\documentclass[11pt]{article}
\usepackage{amssymb, amsmath, amsthm, tikz, xspace,subfigure, graphicx}
\usepackage[left=1in,top=1in,right=1in]{geometry}
\usepackage{stmaryrd}
\newcommand{\eval}[2]{\llbracket #1 \rrbracket_{#2}}
\usetikzlibrary{positioning, shapes.geometric, calc, intersections}
\tikzstyle{vertex}=[circle,fill=black,inner sep=2pt]
\tikzstyle{vertrect}=[draw,rectangle,inner sep=2pt]
\tikzstyle{vertdia}=[draw,diamond,inner sep=2pt]

\theoremstyle{plain}
      \newtheorem{theorem}{Theorem}[section]
      \newtheorem{lemma}[theorem]{Lemma}

      \newtheorem{corollary}[theorem]{Corollary}
      
      \newtheorem{conjecture}[theorem]{Conjecture}
\theoremstyle{definition}
      \newtheorem{definition}[theorem]{Definition}
      \newtheorem{remark}[theorem]{Remark}

\newcommand{\of}[1]{\left( #1 \right)}
\newcommand{\df}{\stackrel{\rm def}{=}}
\newcommand{\Hom}{\text{\rm Hom}}
\newcommand{\function}[2]{:#1 \longrightarrow #2}
\newcommand{\rn}{\boldsymbol}
\newcommand{\prob}[1]{ {\bf P}\! \left[ #1 \right] }
\newcommand{\gkn}{{\cal G}_k(n)}
\newcommand{\gkni}{{\cal G}_k(n_i)}

\def\twr{\mbox{\rm twr}}

\title{Polynomial to exponential transition in Ramsey theory}

\author{Dhruv Mubayi\thanks{Department of Mathematics, Statistics, and Computer Science, University of Illinois, Chicago, IL, 60607 USA.  Research partially supported by NSF grants DMS-1300138 and 1763317. Email: {\tt mubayi@uic.edu}}\and Alexander Razborov\thanks{Departments of Mathematics and Computer Science, University of Chicago, IL, 60637, USA, {\tt razborov@math.uchicago.edu} and Steklov Mathematical Institute, Moscow, 119991, Russia, {\tt razborov@mi.ras.ru}.}}

\begin{document}

\maketitle

\begin{abstract}
Given $s \ge k\ge 3$, let $h^{(k)}(s)$ be the minimum $t$ such that there exist arbitrarily large $k$-uniform hypergraphs $H$ whose independence number is at most polylogarithmic in the number of vertices and in which every $s$ vertices span at most $t$ edges.   Erd\H os and Hajnal conjectured (1972) that $h^{(k)}(s)$ can be calculated  precisely using a recursive formula and Erd\H os offered \$500 for a proof of this. For $k=3$ this has been settled for many values of $s$ including powers of three but it was not known for any $k\geq 4$ and  $s\geq k+2$.

 Here we  settle the conjecture for all $s \ge k \ge 4$.
    We also answer a question of Bhat and R\"odl by constructing, for each $k \ge 4$, a quasirandom sequence of $k$-uniform hypergraphs with positive density and upper density at most $k!/(k^k-k)$. This result is sharp.
\end{abstract}

{\bf MSC classification codes}: 05D10 (primary), 05C35, 05C55, 05C65 (secondary)

\section{Introduction}

Write $K^{(k)}_N$ for the complete $k$-uniform hypergraph (henceforth $k$-graph) on $N$ vertices.
The \emph{Ramsey number} $r_k(s,n)$ is the minimum $N$ such that every red/blue coloring of the edges of $K^{(k)}_N$  contains a monochromatic red copy of $K_s^{(k)}$ or a monochromatic blue copy of $K^{(k)}_n$.
In order to shed more light on the growth rate of these classical Ramsey numbers,  Erd\H os and Hajnal~\cite{EH72} in 1972 considered the  following more general parameter.

\begin{definition} For integers $2\le k < s <n$ and $2 \le t \le {s \choose k}$, let $r_k(s,t;n)$ be the minimum $N$ such that every red/blue coloring of  the edges of $K^{(k)}_N$ contains a monochromatic blue copy of $K_n^{(k)}$ or has a set of $s$ vertices which contains at least $t$ red edges.
\end{definition}
Note that $r_k(s,{s \choose k};n)= r_k(s,n)$ so $r_k(s,t; n)$ includes classical Ramsey numbers.
 In addition, the case
$(k,s,t,n)=(k, k+1, k+1, k+1)$ was investigated in relation to the Erd\H os-Szekeres theorem and Ramsey numbers of ordered tight paths as well as to high dimensional
tournaments by several researchers~\cite{DLR, EM, FPSS, LM, MSW,MS}; the very special case
$(3, 4, 3, n)$ has  connections to quasirandom hypergraph constructions~\cite{BR, KNRS, LM1, LM2}.

The main conjecture of Erd\H os and Hajnal \cite{EH72}  for $r_k(s,t;n)$ is that, as $t$ grows from $1$ to ${s\choose k}$, there is a well-defined value $t_1=h_1^{(k)}(s)$ at which $r_k(s,t_1-1;n)$ is polynomial in $n$ while $r_k(s,t_1;n)$ is exponential in a power of $n$, another well-defined value $t_2=h_2^{(k)}(s)$ at which it changes from exponential to double exponential in a power of $n$ and so on,  and finally a well-defined value $t_{k-2}=h_{k-2}^{(k)}(s)<{s \choose k}$
at which it changes from $\twr_{k-2}$ to $\twr_{k-1}$ in a power of $n$. They were not able to offer a conjecture as to what $h_i^{(k)}(s)$ is in general, except when $i=1$ or when $s=k+1$.

The problem of determining $r_k(k+1, t;n)$ for $t=2$ and $t=3$ has essentially been solved.
For general $t$, the methods of Erd\H os and Rado \cite{ER} show that there exists $c=c(k,t)>0$ such that $r_k(k+1, t; n) \leq \twr_{t-1}(n^c)$ for $3 \leq t\leq k$.  Erd\H os and Hajnal conjectured  that this gives the correct tower growth rate for $ r_k(k+1, t; n)$. When $k\geq 6$, the first author and Suk~\cite{MS16} settled their conjecture in almost all cases in a strong form.

Perhaps the main open problem about $r_k(s, t; n)$ posed by Erd\H os and Hajnal~\cite{EH72} was to determine the value of $t_1=h_1^{(k)}(s)$; namely the value of $t$ at which $r_k(s,t; n)$ transitions from polynomial to super polynomial growth. This is the problem we address in this paper.  The following function plays an important role.

\begin{definition} \label{recursive} Given  positive integers $s,k$, call a partition $s_1+\cdots +s_k=s$ {\em nontrivial} if $0\le s_i < s$ for each $i$. For $0 \le s<k$, let $g_k(s)=0$ and for $s \ge k\ge 3$, let $g_k(s)$ be the maximum of
	$$\sum_{i=1}^k g_k(s_i) + \prod_{i=1}^k s_i$$
	where the maximum is taken over all nontrivial partitions $s_1+\cdots +s_k=s$.
\end{definition}
	We will interpret  $g_k(s)$ as the maximum number of edges in the $s$-vertex $k$-graph obtained by first partitioning $s$ vertices into $k$  parts, taking all edges that intersect all parts, and then recursing this construction within each part.
	Erd\H os and Hajnal commented without proof that it is easy to see that $g_k(s)$ is achieved by taking a partition that is as equitable as possible. We will prove this in the Appendix, and also prove an asymptotic version of this fact later (see \eqref{eq:left_end}). As an easy exercise, this implies that
\begin{equation} \label{eq:asymptotic}
g_k(s)= (1+o(1))\frac{k!}{k^k-k}{s \choose k}\qquad \hbox{$k$ is fixed, $s \to \infty$}.
\end{equation}
		 Erd\H os and Hajnal proved that $r_k(s, g_k(s); n)$ is polynomial in $n$ for all fixed  $s>k \ge 3$. In other words, they showed that every $N$-vertex $k$-graph ($k\ge 3$ fixed) in which every $s$-set spans at most $g_k(s)-1$ edges has independence number at least $N^{\epsilon}$	where $\epsilon>0$ depends only on $s,k$. Therefore
	 $$h_1^{(k)}(s)\ge g_k(s)+1.$$
	
	 They conjectured the following for which Erd\H os later offered \$500 (see~\cite{NR} page 21  and~\cite{Chung} Problem (85)).

\begin{conjecture} [Erd\H os-Hajnal] \label{i=1} Fix $s \ge k \ge 3$. Then $h_1^{(k)}(s)=g_k(s)+1$, or equivalently,
	$r_k(s, g_k(s)+1; n)$ is at least exponential in a power of $n$.
\end{conjecture}
 For $k=3$, Erd\H os and Hajnal~\cite{EH72} proved that Conjecture~\ref{i=1} follows from the following conjecture.

 \begin{conjecture} [Erd\H os-Hajnal] \label{ehrainbow} For every edge-coloring of the complete graph with vertex set $[n]$ by three colors I, II, III,  the  number of triangles  $\{a,b,c\}$ with $a<b<c$ for which $\{a,b\}$ has color I, $\{b,c\}$ has color II, and $\{a,c\}$ has color III is at most $g_3(n)$.
\end{conjecture}
  Conlon-Fox-Sudakov~\cite{CFS} connected Conjecture~\ref{ehrainbow} to  the maximum number $T(s)$ of directed triangles in an $s$-vertex tournament (It is worth noting that the hypergraphs in Conjecture~\ref{ehrainbow} were also considered in~\cite{RRS} due to their connection to hypergraph Tur\'an theory.) They determined $T(s)$ exactly and observed that this also settles Conjecture~\ref{ehrainbow}  for many values of $s$ including powers of 3. Consequently, their approach gave a solution to Conjecture~\ref{i=1} when $k=3$ and these $s$-values; they also proved that $h_1^{(3)}(s)= (1/4){s \choose 3} + O(s \log s)$.  However, their method  using $T(s)$ does not apply for any $k>3$ as it does not capture the recursive structure from Definition~\ref{recursive} needed to prove Conjecture~\ref{i=1}.  Indeed, the set of extremal configurations for $T(s)$  consists of {\em all} (out-)regular tournaments; the recursive construction is just one (and unnecessarily complicated) example in this class. Thus Conjecture~\ref{i=1} was known only when $s=k+1$ and when $k=3$ and $s$ is as described above. In fact, Erd\H os and Hajnal stated in~\cite{EH72} that they were much less certain about Conjecture~\ref{i=1} when $k \ge 4$ than when $k=3$.

 In this paper we prove  Conjecture~\ref{i=1} for all $k\geq 4$.

 \begin{theorem} \label{k}
 	 $h_1^{(k)}(s)= g_k(s)+1$ for all $s \ge k\geq 4$.
 \end{theorem}

Our method also answers a question posed by Bhat and R\"odl~\cite{BR} about quasirandom sequences.
The density of a $k$-graph $H=(V,E)$ is $d(H) = |E|/{|V| \choose k}$. Let ${\cal H} = \{H_n\}_{n=1}^{\infty}$ be a sequence of $k$-graphs with $H_n=(V_n, E_n)$ such that $|V_n| \to \infty$ as $n \to \infty$. Define the density $d({\cal H})$ of ${\cal H}$ as $d({\cal H}) = \lim_{n \to \infty} d(H_n)$ (we only consider sequences where the limit exists) and the upper density
$$\overline{d}({\cal H}) \df \lim_{s \to \infty} \max_n \max_{S \in {V_n \choose s}} d(H_n[S])$$
(note that for any fixed $s$, $H_n[S]$ can take only finitely many values, up to isomorphism).
One can show by a simple averaging argument that $\overline{d}({\cal H})$ exists.

\begin{definition} A $k$-graph sequence $\{H_n\}_{n=1}^{\infty}$ is $\rho$-quasirandom
	 if for every $\epsilon>0$ there exists $n_0$ such that for $n>n_0$, every  $W \subset V = V(H_n)$ with $|W|\ge \epsilon |V|$ satisfies $d(H_n[W])\in [\rho(1-\epsilon),\ \rho(1+\epsilon)]$.
	\end{definition}

An important result of Erd\H os~\cite{E64} states that every $k$-graph sequence with positive density contains arbitrarily large complete $k$-partite subgraphs and hence has upper density at least $k!/k^k$ (the case $k=2$ was done earlier by K\" ov\' ari-S\'os-Tur\'an~\cite{KST} and by Erd\H os (see~\cite{E64})); the value $k!/k^k$ cannot be increased as shown by complete $k$-partite $k$-graphs. This is a fundamental tool for hypergraph problems, and shows that every $\rho \in (0, k!/k^k)$ is a ``jump" for $k$-graphs (see~\cite{FR} for background on jumps).

 Bhat and R\"odl~\cite{BR} improved this result of Erd\H os in the quasirandom setting: they showed that for each $k \ge 3$ and $\rho>0$, every $\rho$-quasirandom $k$-graph sequence has upper density at least  $k!/(k^k-k)$, thus showing that every $\rho \in (0, k!/(k^k-k)$) is a jump in this setting.
It is well-known that $k!/(k^k-k)$ cannot be increased   for $k=3$ (the simplest example is to take the 3-graph of all cyclic triples in a random tournament) and Bhat and R\"odl asked whether the same is true for $k\ge 4$. We answer this positively, showing that the result in~\cite{BR} is sharp for all $k \ge 3$.

\begin{theorem} \label{Rodl}
	For each $k \ge 4$, there exists $\rho>0$ and a $\rho$-quasirandom $k$-graph sequence with  upper density $k!/(k^k-k)$.
	\end{theorem}
We note that our proof of Theorem~\ref{Rodl} yields $\rho = k^{-\Omega(k^2)}$ which is much smaller than $k!/(k^k-k)$ and it remains open to prove the theorem with $\rho = k!/(k^k-k)$ (for $k=3$ this is true).

\section{Reduction to inducibility}

As mentioned in the introduction, Erd\H os and Hajnal showed that the $k=3$ case of Conjecture~\ref{i=1} follows from Conjecture~\ref{ehrainbow}, which asks for the maximum number of rainbow colored triangles (with some additional properties) in an edge-colored ordered graph.  This is an example of a question about the inducibility of colored, directed structures.
 In fact, Erd\H os and Hajnal  observed that Conjecture~\ref{ehrainbow} could be replaced by another slightly different question about inducibility (where we use only two colors and count certain 2-colored triangles) and, as mentioned earlier,  Conlon, Fox and Sudakov~\cite{CFS} considered yet another inducibility problem, namely the determination of $T(s)$.

 Our approach to Conjecture~\ref{i=1} is to formulate a novel question about the inducibility of colored directed structures and solve it exactly. It is perhaps interesting that the
 ``universal" character of the structure we consider below allows us to get around many technical difficulties plaguing previous research on inducibility.
\medskip

\begin{theorem} {\bf (Main Result)} \label{main}
Let $s \ge k\geq 4$ and $R$ be an arbitrary $k$-vertex tournament whose edges are colored with the ${k \choose 2}$ distinct colors from ${[k] \choose 2}$.
	Then the number of copies of $R$ in any $s$-vertex tournament whose edges are colored from ${[k]\choose 2}$ is at most $g_k(s)$.
\end{theorem}

We immediately get Theorem~\ref{k} as a consequence.

\medskip

{\bf Proof of Theorem~\ref{k}.} Fix $s \ge k \geq 4$.  We are to show that $h_1^{(k)}(s) \le g_k(s)+1$. In other words:  there exists $C=C(k)>0$ and, for all $N>k$,  an $N$-vertex $k$-graph $H$ with $\alpha(H) \le C \log N$ such that every $s$ vertices of $H$ span at most $g_k(s)$ edges.
Fix a $k$-vertex tournament $R$ whose edges are colored with ${k \choose 2}$ distinct colors. Next consider the random $N$-vertex tournament $T=T_N$ whose edges are randomly colored with the same ${k \choose 2}$ colors; thus, each pair gets a particular orientation and color with probability $p=1/((k-1)k)$.   Now form the $k$-graph $H=H(T)=(V,E)$ with  $V=V(T)$ and  $E=\{K \subset V: H[K]\cong R\}$. In other words, the edges of $H$ correspond to  copies of $R$. By Theorem~\ref{main}, every $s$ vertices of $H$ span at most $g_k(s)$ edges. On the other hand, the probability that a given $k$-set of vertices in $H$ induces a copy of $R$ is  $k!p^{{k\choose 2}}>0$. Hence the expected number of $t$-sets in $H$ that are independent is at most ${N \choose t}2^{-\epsilon t^2}$ for appropriate $\epsilon=\epsilon(k)>0$.  Indeed, given any $t$-set $A$, pick up in it $\ell\geq\Omega(t^2/k^2)$ $k$-subsets $B_1,\ldots,B_\ell$ such that $|B_i\cap B_j|\leq 1$ whenever $i\neq j$, and notice that the events ``$B_i$ spans a copy of $R$'' are mutually independent. This expectation is less than one as long as $t>C \log N$ and $C=C(k)$ is sufficiently large. \qed
\medskip

\begin{remark}
For the remaining case $k=3$, we believe that the Erd\H os-Hajnal conjecture still holds
but it may require new techniques and ideas: many crucial calculations in this paper completely fall apart.
\end{remark}

\begin{remark}
As mentioned above, there is nothing specific about the kind of combinatorial structures we are considering here, and Theorem \ref{k} is implied by results analogous to Theorem \ref{main} for arbitrary structures.
For example,~\cite{BHLF} gives Theorem \ref{k} for $k=5,\ s=5^t$, \cite{Y} gives it for all $k$ sufficiently large and
$s\leq 2^{\sqrt k}$, the well-known Pippenger-Golumbic conjecture~\cite{PG} about the inducibility of $C_k$ would imply it
for $k\ge 5,\ s=k^t$, and the conjecture from \cite{S} about $\vec C_4$ would imply it for $k=4,\ s=4^t$. See
~\cite{B, BS, EL,FHL, HHN, erdos, H,  KNV} for results about inducibility for other structures.
\end{remark}

Next we show how our approach also answers the question of Bhat and R\"odl about $\rho$-quasirandom hypergraph sequences. It is convenient to use the following theorem from~\cite{CHPS} which is a hypergraph generalization of  the Chung-Graham-Wilson characterization of graph quasirandomness. In what follows $M_k$ is a specific linear $k$-graph with $v=k 2^{k-1}$ vertices and $e=2^k$ edges (see~\cite{CHPS} or~\cite{LM1} for the precise definition); in particular $M_2=C_4$.
We write the result from~\cite{CHPS} in the language of hypergraph sequences.

\begin{theorem} [Conlon-Han-Person-Schacht~\cite{CHPS}] \label{chps} Fix $k \ge 2$, $0<\rho<1$, and a sequence of $k$-graphs ${\cal H} = \{H_n\}_{n=1}^{\infty}$ of density $\rho$ each with $H_n=(V_n, E_n)$ and $|V_n| \to \infty$ as $n \to \infty$. Then ${\cal H}$ is $\rho$-quasirandom iff the number of (labeled) copies of $M_k$ in $H_n$ is
	$|V_n|^{k 2^{k-1}} \rho^{2^k} (1+o(1))$  as $n \rightarrow \infty$.
		\end{theorem}

\medskip

{\bf Proof of Theorem~\ref{Rodl}.}  We use the proof of Theorem~\ref{k} above to construct the desired sequence. Using the notation there, for each $k \ge 4$,  let $\rho=k!((k-1)k)^{-{k \choose 2}}$ be the probability that a $k$-set induces a copy of $R$. For $n \ge 1$, let $\epsilon_n= 1/n$. Standard probabilistic arguments  together with the construction of $H$ in Theorem~\ref{k} imply that there exists a $k$-graph $H_n=(V_n,E_n)$ whose edge set comprises copies of $R$ such that $|V_n| \rightarrow \infty$ and the number of copies of
$M_k$ in $H_n$ is $|V_n|^{k 2^{k-1}} \rho^{2^k} (1\pm \epsilon_n)$. Indeed, since $M_k$ is linear (meaning that every two edges of $M_k$ share at most one vertex) the expected number of (labeled) copies of $M_k$ in $H_n$ is $|V_n|^{k 2^{k-1}} \rho^{2^k}$ and Chebyshev's inequality implies that there is an $H_n$ where the number of copies of
$M_k$ in $H_n$ is $|V_n|^{k 2^{k-1}} \rho^{2^k} (1\pm \epsilon_n)$.   Now let ${\cal H} = \{H_n\}_{n=1}^{\infty}$. We have just shown that  the number of  copies of $M_k$ in $H_n$ is $|V_n|^{k 2^{k-1}} \rho^{2^k} (1+o(1))$ so Theorem~\ref{chps} implies that
${\cal H}$ is  $\rho$-quasirandom. On the other hand, for each $s,n>0$, $k \ge 4$ and $S \in {V_n \choose s}$ we have $d(H_n[S]) \le g_k(s)/{s \choose k}$ by Theorem~\ref{main}. Consequently, $\overline{d}({\cal H})\le \lim_{s \to \infty} g_k(s)/{s \choose k} \le k!/(k^k-k)$.  \qed

\section{Proof of asymptotic result} \label{sec:asymptotic}

Recall that the {\em inducibility $i(R)$} is defined as
$$
i(R) \df \lim_{s\to\infty} \max_{|V(H)|=s} \frac{i(R;H)}{{s \choose k}},
$$
where $i(R;H)$ is the number of copies of $R$ in an $s$-vertex ${[k]\choose 2}$-colored tournament $H$. In this section we prove the following.

\begin{theorem} \label{asymptotic} Let $k\geq 4$ and $R$ be an arbitrary $k$-vertex tournament whose edges are colored with the ${k \choose 2}$ distinct colors from ${[k] \choose 2}$. Then
$$
i(R) = \frac{k!}{k^k-k}
$$ {\rm (}which, by \eqref{eq:asymptotic} is equal to $\lim_{s\to\infty}\frac{g_s(k)}{{s\choose k}}${\rm )}.
\end{theorem}

The proof of Theorem~\ref{asymptotic} is much cleaner than the proof of our main result $\max_{|V(H)|=s} i(R;H) = g_k(s)$ presented in Section \ref{sec:exact} since it avoids dealing with unnecessary details about the number of vertices. It also gives the reader the overall structure of our argument.  Moreover, as we will show in Corollary~\ref{kpowers}, the asymptotic result in Theorem~\ref{asymptotic} immediately implies an exact result whenever $s$ is a power of $k$.

To make our argument both clean and rigorous, we use the language of Flag Algebras. But since in order to prove Theorem \ref{main} we will have to ``discretize'' it anyway (so Theorem \ref{asymptotic} is sort of a warm-up), we skip the traditional crash course in Flag Algebras and assume a certain degree of familiarity with the method. The reader interested only in the end result can safely proceed to Section \ref{sec:exact} (or, if willing to believe that all this can be made completely rigorous, follow the proof on the intuitive level).

\medskip
{\bf Proof of Theorem~\ref{asymptotic}.} Let $T_k$ be the theory \cite[\S 2]{flag} of ${[k] \choose 2}$-colorings
of edges of a complete graph, and let $T$ be the disjoint union of $T_k$ and the theory  $T_{\text{Tournament}}$ of tournaments. Let $R\in \mathcal M_k[T]$ be any model with $V(R) = [k]$ such that its restriction to $T_k$ is the canonical (that is, the edge $(i,j)$ is colored with the color $\{i,j\}$) model, and let $\Gamma$ be the underlying tournament. As always, we denote by 1 the (only) type \cite[\S 2.1]{flag} of size 1.

For a color $c\in {[k]\choose 2}$, there are two 1-flags in $\mathcal F_2^1$ colored
by $c$: $\alpha_c$ (in which the distinguished vertex is the tail) and
$\beta_c$ (distinguished = head). Let
$$
S_i\df\sum_{j\in
N_\Gamma^+(i)}\alpha_{\{i,j\}}+ \sum_{j\in
N_\Gamma^-(i)}\beta_{\{i,j\}}
$$
(this is an element of $\mathcal A^1$) and note that $\sum_i S_i=1$. Define also
\begin{equation} \label{eq:p_i}
P_i\df\prod_{j\in N_\Gamma^+(i)}\alpha_{\{i,j\}}\cdot \prod_{j\in
N_\Gamma^-(i)}\beta_{\{i,j\}}.
\end{equation}

Let us now fix $\phi\in\Hom^+(\mathcal A^0[T],\mathbb R)$ (see \cite[Definition 5]{flag}) maximizing $\phi(R)$ \cite[\S 4.1]{flag}
and let
$$
a_k\df \frac{k!}{k^k-k}\ \of{=\frac{(k-1)!}{k^{k-1}-1}}.
$$
Our goal is to prove that
\begin{equation} \label{goal}
\phi(R)\leq a_k,
\end{equation}
and we can assume w.l.o.g. that $\phi(R)\geq a_k$.

Let $\rn{\phi^1}$ be the distribution  over $\Hom^+(\mathcal A^1[T], \mathbb R)$ rooted at $\phi$ \cite[Definition 10]{flag}, and let $\mathcal S^1(\phi)$ be the support of this distribution. Combinatorially, $\rn{\phi^1}$ should be thought of as a uniform distribution over vertices (except that we do have any such thing as a “vertex” here). Let us study an individual element $\phi^1\in \mathcal S^1(\phi)$.

Assume for simplicity that $\phi^1(S_1)\geq \phi^1(S_2)\geq\cdots\geq\phi^1( S_k)$; our goal is to bound $\phi^1(S_2)$ from above (the trivial bound is $1/2$). More specifically, note first (recall that $k\geq 4$) that
$$
\frac{(k-1)^{k-1}}{k^{k-1}-1} \geq \of{1-\frac 1{k}}^{k-1} \geq e^{-1}\geq 2^{2-k}.
$$
Hence the equation
\begin{equation} \label{roots}
z^{k-1} + (1-z)^{k-1} = \frac{(k-1)^{k-1}}{k^{k-1}-1}
\end{equation}
has two roots in the interval $(0,1)$; let $z_k\in (0,1/2)$ be the smallest one. We claim that $\phi^1(S_2)\leq z_k$.

Let  $\mu^1_k(R)\in \mathcal A^1$ be the sum of all $k$ possible 1-flags that can be obtained from $R$ \cite[\S 4.3]{flag}. Then one consequence of the
extremality of $\phi$ (and the fact that $\phi^1\in \mathcal S^1(\phi)$) is that $\phi^1(\mu^1_k(R))=\phi(R)\geq a_k$ \cite[Theorem 4.3]{flag}.

On the other hand, by the AMGM inequality we have
\begin{equation} \label{mu_bound}
\mu^1_k(R) \leq (k-1)!\sum_i P_i \leq \frac{(k-1)!}{(k-1)^{k-1}}\sum_i S_i^{k-1},
\end{equation}
where the partial pre-order $\leq$ on $\mathcal A^1$ simply means \cite[Definition 6]{flag} that the inequality holds upon being
evaluated by an arbitrary element of $\Hom^+(\mathcal A^1[T],
\mathbb R)$.
Comparing the two,
$$
\sum_i \phi^1(S_i)^{k-1} \geq \frac{(k-1)^{k-1}}{k^{k-1}-1}.
$$
But under the condition $S_2=z\ (\leq 1/2)$, the left-hand side is clearly maximized when $S_1=1-z$ and $S_3=\ldots=S_k=0$. This gives us the claim.

We now have a measurable partition $S^1(\phi) = V_1\stackrel .\cup\ldots \stackrel .\cup
V_k$ according to $\text{arg max}_i\phi^1(S_i)$ (we resolve conflicts arbitrarily), and we want to incorporate  it into our language explicitly.
Let $T^+$ be the extension of
$T$ with vertex coloring $\chi$ in $k$ colors. We let
$p_i\in\mathcal M_1[T^+]$ be the one-element model in which the only vertex
is colored by $i$, and let $(i)$ be the corresponding type. Let $I: T\leadsto T^+$ be the interpretation \cite[\S 2.3.3]{flag} erasing vertex coloring. We want to extend $\phi$ to an element $\phi^+\in\Hom^+(\mathcal A^0[T^+],\mathbb R)$ that respects the partition $V_1\stackrel .\cup\ldots \stackrel .\cup
V_k$ (we will actually need only its property $\phi^1(S_{i})\leq z_k$ for $\phi^1\not\in V_i$). Formally, we claim the existence of $\phi^+$ with the
following two properties:
\begin{enumerate}
\item $\phi = \phi^+\circ \pi^I$ (for the definition of algebra homomorphisms $\pi^I, \pi^{I,\sigma},\pi^{\sigma,\eta}$ etc. see \cite[\S 2.3]{flag});

\item For any $\psi\in\mathcal S^{(i)}(\phi^+)$ and any $i'\neq i$, $\psi\of{\pi^{I,(i)}(S_{i'})}\leq z_k$.
\end{enumerate}
Combinatorially, the existence of such an extension is entirely obvious, and the simplest way to give a rigorous definition in the language of Flag Algebras is by an explicit formula. Namely, for a type $\sigma$ of the theory $T^+$ that has size $k$, we first define the ``labelled density'' $\phi^+(\langle\sigma\rangle)$ as
\begin{equation} \label{eq:extension}
\phi^+(\langle\sigma\rangle)\df \phi(\langle I(\sigma)\rangle)\cdot\prob{\bigwedge_{i=1}^k\of{\rn{\phi^{I(\sigma)}}\circ
\pi^{I(\sigma),i}\in V_i}}.
\end{equation}
Then we let
$$
\phi^+(\sigma)\df (S_k: \text{Aut}(\sigma))\phi^+(\langle\sigma \rangle).
$$
It is straightforward to check that so defined $\phi^+$ is an element of
$\Hom^+(\mathcal A^0[T^+], \mathbb R)$ that satisfies properties 1) and 2) above.

From now on we will often omit from the notation operators $\pi^I$ and $\pi^{I,\sigma}$ (as well as $\phi^+$); thus, the algebra $\mathcal A^{I(\sigma)}[T]$ is identified with its
image under $\pi^{I,\sigma}$ in $\mathcal A^\sigma[T^+]$. When $\sigma$ has to be specified, we write $f^\sigma$ for the image of $f\in \mathcal A^{I(\sigma)}[T]$ in
$\mathcal A^\sigma[T^+]$; we will be primarily interested in the case when $\sigma$ has size 1, i.e. $\sigma=(i)$ for some $i\in [k]$. Thus, property  2) above simplifies to $\psi\of{S_{i'}^{(i)}}\leq z_k$ for any $i'\neq i$ and $\psi\in\mathcal S^{(i)}(\phi^+)$ etc.

For $j\neq i$, let $P_{ij}\in\mathcal A^1[T]$ be the product $P_i$ with either $\alpha_{i,j}$ or
$\beta_{i,j}$ removed. Then the AMGM inequality implies the bound
\begin{equation}\label{bound_on_hat}
  \psi\of{P_{i'j}^{(i)}}\leq \of{\frac{z_k}{k-2}}^{k-2}\ \of{i'\neq i,\ \psi\in\mathcal S^{(i)}(\phi^+)}.
\end{equation}

Now, $R$ splits in $T^+$ as follows:
$$
\pi^I(R)=R_m + R_g + R_b
$$
("m, g, b" stand for ``monochromatic'', ``good'' and ``bad'', respectively),
where $R_m$ is the sum of $m$ models in $I^{-1}(R)$ in which all vertices are colored in the same color, $R_g$ is the model with $\chi=\text{id}$ and $R_b$ is the sum of all remaining models. We will estimate these three terms (evaluated by $\phi^+$) separately.

\smallskip
The bound on $R_m$ (that, combinatorially, is the density of monochromatic copies of $R$) is obtained by exploiting the extremality of $\phi$ one more time:
\begin{equation} \label{phi_plus}
\phi^+(R_m) \leq \phi(R) \cdot \sum_i \phi^+(p_i)^k.
\end{equation}

To make this rigorous, whenever $\phi^+(p_i)>0$, we form the restriction $\pi^{p_i}\function{\mathcal A^0[T^+]}{\mathcal A^0_{p_i}[T^+]}$, where $\mathcal A^0_{p_i}$ is the localization of $\mathcal A^0$ by the element $p_i$ (see \cite[\S 2.3.2]{flag} with $\sigma=0$) that combinatorially corresponds to the restriction on $V_i$. Then $\phi^+\circ\pi^{p_i}\circ \pi^I\in\Hom^+(\mathcal A^0[T],\mathbb R)$, hence the extremality of $\phi$ implies that $(\phi^+\circ\pi^{p_i}\circ \pi^I)(R)\leq \phi(R)$. On the other hand, by unrolling definitions we see that $(\phi^+\circ\pi^{p_i}\circ \pi^I)(R)=\frac{\phi^+(R_{m,i})}{\phi^+(p_i)^k}$, where $R_{m,i}$ is the model with $\chi\equiv i$. Multiplying by $p_i^k$ and summing over all $i$ gives us \eqref{phi_plus}.

\smallskip
For bounding $R_g$, we let $\Delta$ be the sum of all rainbow (i.e., with bijective $\chi$) models in $\mathcal M_k[T^+]$ {\em different} from $R_g$. Then we clearly have
\begin{equation} \label{eq:r_decomposition}
R_g = k!\prod_{i=1}^k p_i-\Delta,
\end{equation}
and we need to bound $\Delta$ from below. We will do it in terms of the element $\delta\in\mathcal A^0_2[T^+]$ which is the sum of all models $M$ with the set of vertices $\{u,v\}$ that are {\em transversal} ($\chi(u)\neq \chi(v)$) and are either {\em miscolored} (the edge color of $(u,v)$ is different from $\{\chi(u), \chi(v)\}$) or {\em disoriented} (the orientation of $(u,v)$ is different from the orientation of $(\chi(u),\chi(v))$ in $\Gamma$), or both. Let $\delta_{ij}$ be the contribution to $\delta$ made by those models in $\mathcal M_2[T^+]$ for which $\{\chi(u),\chi(v)\} = \{i,j\}$. Then
$$
\delta = \sum_{1\leq i<j\leq k} \delta_{ij}.
$$
Now, if we extend any model in $\delta_{ij}$ to an arbitrary rainbow model in $\mathcal M_k[T^+]$, we will actually get a model in $\Delta$. This implies
$$
\Delta\geq \frac{k!}{2}\of{\prod_{\nu\neq i,j}p_\nu} \delta_{ij}
$$
for any $i\neq j$ (the factor 2 in the enumerator accounts for the symmetry interchanging $i$ and $j$). Multiplying this by $p_ip_j$ and summing up over all such pairs, we arrive at
\begin{equation} \label{Delta_estimate}
\of{\sum_{1\leq i<j\leq k}p_ip_j}\Delta \geq \frac{k!}{2}\of{\prod_i p_i}\delta.
\end{equation}

\medskip
At this point we have to take care of the case when all but one of the $p_i$s are equal to 0. Assuming that, say, $p_1=1$, we know that $\phi^1(S_2)\leq z_k$ for {\em any} $\phi^1\in S^1(\phi)$. Let $R^{1,i}\in\mathcal F^1[T]$ be obtained from $R$ by placing the distinguished vertex into $i$ (so that $\mu^1_k(R)=\sum_i R^{1,i}$). Then the local version of \eqref{mu_bound} gives us
\begin{equation} \label{eq:r12}
R^{1,2}\leq \frac{(k-1)!}{(k-1)^{k-1}} S_2^{k-1} \leq \frac{(k-1)!}{(k-1)^{k-1}} z_k^{k-1}.
\end{equation}
On the other hand, $R=k\eval{R^{1,2}}{1}$ (see \cite[\S 2.2]{flag} for the averaging operator $\eval{\cdot}{\sigma}$). This gives us
\begin{equation} \label{eq:r}
\phi(R) \leq \frac{k!}{(k-1)^{k-1}}z_k^{k-1} \leq a_k k\of{1+\frac 1{k-1}}^{k-1}z_k^{k-1} \leq a_k,
\end{equation}
where for $k=4$ the last inequality follows from the bound $z_4\leq 0.3$, and when $k\geq 5$ it suffices to apply the trivial bound $z_k\leq 1/2$. This completes the proof of \eqref{goal} when $p_1=1$.

Thus we can and will assume that $p_i<1$ for all $i$ and hence we can divide \eqref{Delta_estimate} by $\sum_{1\leq i<j\leq k}p_ip_j$. Comparing the result with \eqref{eq:r_decomposition}, we arrive at our second estimate
\begin{equation} \label{eq:second_estimate}
R_g\leq k!\cdot \prod_i p_i\of{1-\frac{\delta}{2\sum_{1\leq i<j\leq k}p_ip_j}} = k!\cdot \prod_i p_i\of{1-\frac{\delta}{1-\sum_ip_i^2}}.
\end{equation}

Let's now turn to upper bounding $R_b$. We need a few simple remarks first.

Every model $M\in\mathcal M_2[T]$ has a unique embedding $\alpha_M\function MR$; intuitively, this embedding corresponds to the ``intended'' vertex-coloring of $M$. The mapping $\alpha$ can be extended to $\mathcal M_2[T^+]$ simply by letting $\alpha_M\df \alpha_{I(M)}$ i.e. by ignoring the vertex-coloring $\chi$. Then $\delta$ can be described as the sum of all transversal models
in $\mathcal M_2[T^+]$ for which $\chi\neq\alpha_M$. Let now $b^{(i)}$ be the sum of all $(i)$-flags that have the form $(M,v)$, where $M$ appears in $\delta$ and $v\in V(M)$ is miscolored, that is $\chi(v) (=i)\neq\alpha_M(v)$
. This element further splits as
$$
b^{(i)} = \sum_{i'\neq i\atop j} b^{(i)}_{i',j},
$$
where $b^{(i)}_{i',j}$ consists of those $(M,v)$ for which $\alpha_M(v)=i'$
and $\alpha_M(u)=j$ ($V(M)=\{u,v\}$). Going one step further,
$$
b^{(i)}_{i',j} = b^{(i)'}_{i',j} + b^{(i)''}_{i',j},
$$
where $b^{(i)'}_{i',j}$ consists of those $(M, v)$ in which the second vertex $u$ is also miscolored, that is $\chi(u)\neq j$. Note for the record that
\begin{equation} \label{eq:b_vs_delta}
 \delta=\sum_i\eval{\sum_{i'\neq i\atop j}b^{(i)'}_{i'j}+2b^{(i)''}_{i'j}}{(i)}
\end{equation}
(the extra coefficient 2 balances off the coefficient $\frac 12$ that will appear in $\sum \eval{b^{(i)''}_{i'j}}{(i)}$ for those models $M\in\delta$ in which only one vertex is miscolored).

The upper bound on $R_b$ will be actually given in terms of the expression $\sum_i\eval{\sum_{i'\neq i\atop j}(b^{(i)'}_{i'j}+2b^{(i)''}_{i'j})P^{(i)}_{i'j}}{(i)}$ differing from the right-hand side of \eqref{eq:b_vs_delta} only in the extra term $P^{(i)}_{i'j}$. For that we have to bound from below $p^M\of{\sum_i\eval{\sum_{i'\neq i\atop j}(b^{(i)'}_{i'j}+2b^{(i)''}_{i'j})P^{(i)}_{i'j}}{(i)}}$ (see \cite[Definition 7]{flag} for $p^M$), where $M$ is a model of size $k$ appearing in $R_b$. This quantity, however,
has a very clean combinatorial meaning. Namely, let $c(M)$ be the number of ordered pairs $\langle i',j\rangle$ such that $\chi(i')\not\in\{i',\chi(j)\}$, where those pairs for which $\chi(j)=j$ are counted twice. Then we have

\begin{equation} \label{eq:p_m}
p^M\of{\sum_i\eval{\sum_{i'\neq i\atop j}(b^{(i)'}_{i'j}+2b^{(i)''}_{i'j})P^{(i)}_{i'j}}{(i)}} = \frac{c(M)}{k!}.
\end{equation}
The reason is simply that any pair $\langle i',j\rangle$ as described above determines an embedding of either $b^{(i)'}_{i'j}$ or $b^{(i)''}_{i'j}$ into $M$, with an appropriate coefficient. But once it is determined, there is precisely one way of assigning the remaining $(k-2)$ vertices to terms in the product $P_{i'j}$ (which is simply \eqref{eq:p_i} with $i:=i'$ and the term corresponding to $\{i',j\}$ missing).

We claim that $c(M)\geq 2(k-2)$. Indeed, another way to interpret $c(M)$ is as twice the number of {\em unordered} pairs $\{i,j\}$ that are transversal ($\chi(i)\neq \chi(j)$) and in which {\em at least one} of the two vertices is miscolored. Now, if $\chi$ is a (non-identical) permutation then the transversality restriction becomes void. Picking arbitrarily any miscolored $i$ and any $j\neq i$ will already give us $(k-1)$ pairs of the desired form. If, on the other hand, $\chi$ is not a permutation, let $C$ be any non-trivial $\chi$-colored class: $2\leq |C|\leq k-1$ (the latter condition holds since  $\chi\neq\text{const}$). At least $|C|-1$ vertices in this class are miscolored which gives us $\geq (|C|-1)(k-|C|)\geq k-2$ desired pairs.

Thus, by \eqref{eq:p_m}, \eqref{bound_on_hat} and \eqref{eq:b_vs_delta} we conclude that
\begin{eqnarray}
\nonumber R_b &\leq& \frac{k!}{2(k-2)} \sum_i\eval{\sum_{i'\neq i\atop j}(b^{(i)'}_{i'j}+2b^{(i)''}_{i'j})P^{(i)}_{i'j}}{(i)} \leq \frac{k!}{2(k-2)} \of{\frac{z_k}{k-2}}^{k-2} \sum_i\eval{\sum_{i'\neq i\atop j}b^{(i)'}_{i'j}+2b^{(i)''}_{i'j}}{(i)}\\ &=& \delta \frac{k!}{2(k-2)} \of{\frac{z_k}{k-2}}^{k-2}. \label{eq:bound_on_rb}
\end{eqnarray}

Along with \eqref{phi_plus} and \eqref{eq:second_estimate}, this gives us
$$
\phi(R)\leq \phi(R)\sum_i p_i^k + k!\prod_i p_i\of{1-\frac{\delta}{1-\sum_ip_i^2}} +  \frac{\delta k!}{2(k-2)}\of{\frac{z_k}{k-2}}^{k-2}.
$$
But since the case $\sum_i p_i^k=1$ was already treated above, in order to finish the proof of \eqref{goal}, it remains to show that
$$
a_k\of{1-\sum_i p_i^k} \geq k! \prod_i p_i\of{1-\frac{\delta}{1-\sum_ip_i^2}} +  \frac{\delta k!}{2(k-2)}\of{\frac{z_k}{k-2}}^{k-2}
$$
or, cancelling the factorials,
$$
\frac{1-\sum_i p_i^k}{k^k-k} \geq \prod_i p_i\of{1-\frac{\delta}{1-\sum_ip_i^2}} +  \frac{\delta}{2(k-2)}\of{\frac{z_k}{k-2}}^{k-2}.
$$

We now make use of the fact that the right-hand side here is linear in $\delta$, hence it suffices to check
our inequality at the end-points of the interval $\delta\in\left[0, 1-\sum_ip_i^2\right]$.

\smallskip
The left end $\delta=0$ amounts to
\begin{equation} \label{eq:left_end}
(k^k-k)\prod_i p_i + \sum_ip_i^k \leq 1.
\end{equation}
Let $S\df [k]^k\setminus\{(i, \ldots, i): i \in [k]\}$. Since $|S|=k^k-k$, the AMGM inequality implies that
$$\sum_{(i_1, \ldots, i_k) \in S} p_{i_1}\cdots p_{i_k}\ge (k^k-k)\left(\prod_{(i_1, \ldots, i_k) \in S} p_{i_1}\cdots p_{i_k}\right)^{\frac{1}{k^k-k}}=(k^k-k)\prod_{i} p_i.$$
Adding $\sum_ip_i^k$ to both sides gives us the desired inequality since
$$
\sum_{(i_1, \ldots, i_k) \in S} p_{i_1}\cdots p_{i_k} + \sum_ip_i^k = \of{\sum_i p_i}^k=1.
$$

\smallskip
The right end $\delta= 1-\sum_ip_i^2$ leads to
$$
\frac{1-\sum_ip_i^k}{1-\sum_ip_i^2} \geq \frac {(k^k-k)}{2(k-2)}\of{\frac{z_k}{k-2}}^{k-2}.
$$

For $k\geq 5$ we simply use the trivial bound $1-\sum_ip_i^k\geq 1-\sum_ip_i^2$ so that we have to prove
$\frac {k^k-k}{2(k-2)}\of{\frac{z_k}{k-2}}^{k-2}
\leq 1$. For $k=5,6$ we can use numerical bounds\footnote{A simple Maple worksheet verifying this fact, as well as several other facts of similar nature below, can be found at {\tt http://homepages.math.uic.edu/\~{}mubayi/papers/ErdosHajnalmw.pdf} and  {\tt http://people.cs.uchicago.edu/\~{}razborov/files/ErdosHajnal.mw}}
$z_5,z_6\leq 0.3$ and for $k\geq 7$ the trivial bound $z_k\leq 1/2$ suffices.

When $k=4$, we have to do a bit of extra work. Assume w.l.o.g. that $p_1$ is the largest. Then $\sum_ip_i^4\leq p_1^2\cdot \sum_ip_i^2$ and
$$
\frac{1-\sum_i p_i^k}{1-\sum_ip_i^2} \geq \frac{1-p_1^2\sum_i p_i^2}{1-\sum_ip_i^2} \geq \frac 43(1-p_1^2/4)
$$
(the latter inequality follows from $\sum_ip_i^2\geq 1/4$). Since $z_4\leq 0.257$, all that remains to prove is $p_1\leq 0.91$. For that we simply re-use our previous calculation showing that w.l.o.g. we can assume $p_1<1$.
Indeed, \eqref{eq:r12} is still true for the flag\footnote{A slightly better bound will be obtained in Lemma \ref{pi}; the improvement is achieved by selecting the {\em minimal} flag among $(R^{1,2})^{(1)},\ldots,(R^{1,k})^{(1)}$
rather than arbitrary. But we need not be that precise here.} $(R^{1,2})^{(1)}$. That is, under the additional assumption that the distinguished vertex is in $V_1$ we have
$$
(R^{1,2})^{(1)} \leq \frac 29z_4^3,
$$
and for all other $i$ we still have the same bound but with the trivial estimate $S_2\leq 1$:
$$
(R^{1,2})^{(1)} \leq \frac 29\ (2\leq i\leq k).
$$
Now the bound \eqref{eq:r} reads as
$$
\phi(R) \leq \frac 89(z_4^3p_1 +1-p_1).
$$
Since $z_4\leq 0.257$, this is $\leq\frac 2{21}\ (=a_4)$ whenever $p_1\geq 0.91$. Hence we can assume w.l.o.g. that $p_1\leq 0.91$ and, as we already observed, this implies \eqref{main} for $k=4$.

This completes the proof of Theorem \ref{asymptotic}. \qed

\begin{corollary} \label{kpowers}
	$h_1^{(k)}(s)=g_k(s)+1$ for all $k \ge 4$ whenever $s$ is a power of $k$.
	\end{corollary}
\proof Let $H$ be a model of $T$ with $|V(H)|=s$. We have to prove that ${{n\choose k}}i(R;H)\leq g_k(s)$. The plan is clear (and well-known): turn $H$ into an element of $\Hom^+(\mathcal A^0[T],\mathbb R)$ by replacing every $v\in V(H)$ with the infinite recursive construction and then apply Theorem \ref{asymptotic} to it. There are several ways to make this intuition rigorous: we might consider convergent sequences or simply come up with an explicit formula as in \cite[Section 2.3]{erdos}. Let us do it geometrically (cf. \cite[Section 2]{ch}) as this is the most elegant one.

Consider the infinite lexicographic product $\Omega\df H\times R^{\infty}$. Thus, the vertices are infinite sequences ${\sf x}=(x_0,x_1,\ldots,x_n,\ldots)$, where $x_0\in V(H)$ and $x_i\in V(R)\ (i\geq 1)$. The edge coloring and the orientation between ${\sf x}\neq {\sf y}$ are read from the first coordinate $i$ in which $x_i\neq y_i$. Further, $\Omega$ is equipped with the measure that is the product of uniform measures on $V(H),V(R)$ and all this structure turns $\Omega$ into a $T$-on (\cite[Definition 3.2]{theon}). Hence we also have (\cite[Theorem 6.3]{theon}) the corresponding algebra homomorphism $\phi\in \Hom^+(\mathcal A^0[T],\mathbb R)$; its values are computed as obvious integrals over $\Omega$. In particular, $\phi(R)$ is given by the ``expected'' formula
$$
\phi(R) = \frac{(s)_k}{s^k}i(R;H) + \frac{a_k}{s^{k-1}}.
$$
Along with Theorem \ref{asymptotic}, this leads us, after a bit of manipulations, to the
bound
$$
{{s\choose k}}i(R;H) \leq \frac{s^k-s}{k^k-k}.
$$
When $s$ is a power of $k$, the right-hand side here is exactly $g_k(s)$ (by an obvious induction).
\qed		
		
\section{Proof of Theorem~\ref{main}} \label{sec:exact}
Before commencing with the formal proof of Theorem~\ref{main} we state some facts about partitions. Recall that a  partition $n_1+\cdots + n_k=n$ is {\em equitable} if
$|n_i-n_j|\le 1$ for all $i \ne j$.

\begin{definition} Let $p(0,0)=1$ and for $q>t>0$, let $p(q,t)$  be the maximum of  $\prod_{i} q_i$ where $q_1+\cdots +q_t=q$ is a partition of $q$ and each $q_i<q$.
\end{definition}
It is easy to see that this maximum is achieved only by an equitable partition.

The following Lemma was stated by Erd\H os and Hajnal~\cite{EH72}. Since we could not find a proof of this, we will give a proof in the Appendix.
\medskip

\begin{lemma} \label{trivial}
	If $n\ge k\ge 3$, then $g_k(n)$ is achieved by an equitable partition.
\end{lemma}

An immediate consequence of Lemma~\ref{trivial} is that
\begin{equation}\label{small_n}
g_k(n) = p(n,k) \qquad \hbox{for $n \le k(k-1)$.}
\end{equation}
Indeed, for $n\le k(k-1)$ the equitable partition for $n$ has each part of size less than $k$.

The next simple lemma collects some useful facts about $p(n,k)$. Its easy proof is left to the reader.
\medskip

\begin{lemma} \label{basics} Let $k \ge 1$.
	\begin{enumerate}
		\item[a)] $p(n+1,k)-p(n,k) = p(n-\lfloor n/k\rfloor, k-1)$.
		\item[b)] $p(n,k)$ is {\bf strictly} increasing whenever $n\geq k-1$.
		\item[c)] for $n\geq n'\geq 1$,
		$$
		p(n+1,k) +p(n'-1,k) \geq p(n,k) +p(n',k).
		$$
	\end{enumerate}
\end{lemma}

{\bf Proof of Theorem~\ref{main}}. In our proof we try to keep the notation reasonably consistent with Section \ref{sec:asymptotic} although some differences are unavoidable.

Fix $k \ge 4$ and a $k$-vertex tournament $R$ with vertex set $[k]$ and pair $\{i,j\}$ is  colored by $\{i,j\}$ from $C={[k] \choose 2}$. Let  $H$ be  an $n$-vertex tournament with edges colored from $C$ and let $i(R; H)$ be the number of copies of $R$ in $H$. We are to prove that  $i(R; H) \le g_k(n)$ and we will proceed by induction on $n$.

 For a vertex $x$ in $V(H)$ and $i \in [k]$, write $d_i(x)$ for the number of copies of $R$ containing $x$ where $x$ plays the role of vertex $i$ in $R$.
 More formally, $d_i(x)$ is the number of isomorphic embeddings $\phi: R \to H$ such that $\phi(i) = x$.
  Let $d(x)=\sum_i d_i(x)$ be the number of copies of $R$ containing $x$. For $i \in [k]$, let $N_i(x)$ be the set of those $y \in V(H)\setminus \{x\}$ for which there is a copy of $R$ in $H$ containing both $x$ and $y$ in which $x$ plays the role of vertex $i$ in $R$. Due to uniqueness of the colors of $R$ we have $N_j(x) \cap N_{j'}(x) =\emptyset$ for $j \ne j'$. Moreover, $N_i(x)$ also has a (unique) partition $\cup_{j \ne i} N_i^j(x)$ where $N_i^j(x)$ comprises those $y$ such that $x,y$ lie in a copy of $R$ with $x$ playing the role of $i$ and $y$ playing the role of $j$\ \footnote{Thus, in the language of Section \ref{sec:asymptotic}, the flag $S_i$ provides a simple upper bound on the density of $N_i(x)$ while $\alpha_{\{i,j\}}/\beta_{\{i,j\}}$ upper bound $N_i^j(x)$.}. We have
  \begin{equation} \label{eq_dn}
  d(x) = \sum_{i=1}^k d_i(x) \le \sum_{i=1}^k \prod_{j \ne i} |N_i^j(x)| \le \sum_{i=1}^k  p(|N_i(x)|, k-1).
  \end{equation}

We now partition $V(H)$ into $V_1 \cup \cdots \cup V_k$, $n_i=|V_i|$, where
$$V_i = \{x \in V(H):  d_i(x) \ge d_j(x) \hbox{ for all } j \ne i\}
 \footnote{The corresponding definition of $V_i$ in Section \ref{sec:asymptotic} considers $|N_i(x)|$ instead of $d_i(x
	)$.}$$ and subject to this property, minimize $\sum_{i,j} |n_i-n_j|$. Note that $n_i<n$ for all $i$, for if, say, $n_1=n$, then $d_1(x) \ge d_2(x)$ for all $x$, and $\sum d_1(x) = i(R; H) = \sum d_2(x)$ implies that  $d_1(x)= d_2(x)$ for all $x$ so we could move a vertex to $V_2$, contradicting the choice of partition.
\medskip

{\bf Claim 1.} $i(R; H) \le  g_k(n)$ for all $n \le k(k-1)$.
\proof We proceed by induction on $n$; the case $n\le k$ is trivial. Pick a
vertex $x$ in $H$ and suppose that there are $i\ne j$ with both $d_i(x)$ and
$d_j(x)$ positive. Let $N_{i} = N_{i}(x)$ and $m_i=|N_i|$. Then (\ref{eq_dn}) gives
$d(x) \le \sum_i  p(m_i, k-1)$. By Lemma \ref{basics}c),
this is maximized when there exist $1\le i< j\le k$ with $m_i+m_j=n-1$ and, since $d_i(x), d_i(x)>0$, one of $m_i, m_j$ is equal to $k-1$. Consequently,
$$d(x) \le p(m_i, k-1) + p(m_j, k-1) \le 1+p(n-k, k-1).$$
Deleting $x$ we have, by induction, at most $i(R; H-x) \le g_k(n-1)$ copies of $R$ and hence $i(R; H) \le i(R; H-x) + d(x) \le g_k(n-1)+ 1+p(n-k, k-1)$. We
claim that
\begin{equation}\label{induction}
g_k(n-1)+ 1 + p(n-k, k-1) \le g_k(n)
\end{equation}
for $n\leq k(k-1)$. Indeed, applying \eqref{small_n} and Lemma
\ref{basics}a), this is equivalent to
$$
p(n-1-\lfloor (n-1)/k\rfloor, k-1) \geq 1+ p(n-k,k-1)
$$
 which in turn follows
from $n-1-\lfloor (n-1)/k\rfloor\geq \max(k-1,n-k+1)$ by Lemma
\ref{basics}b).

We may now assume  that for each vertex $x$ there is a unique $i$ for which
$d_i(x)>0$ (otherwise apply \eqref{induction}) and this gives a natural
$k$-partition of the vertex set. Moreover, we now also have $i(R; H)\le p(n, k)=g_k(n)$ by (\ref{small_n}) since $n \le k(k-1)$. \qed

Claim 1 concludes the base case and we now proceed to the induction step where we may assume that $n > k(k-1)$.  We may also assume that  $d(x) \ge d_{min}=g_k(n)-g_k(n-1)$ for each vertex $x$ as otherwise we may delete $x$ and apply induction.

The next part of the argument (up to the inequality \eqref{final}) closely parallels the one given in Section \ref{sec:asymptotic} but we give it here anyway for the sake of completeness.

Partition the copies of $R$ in  $H$  as $H_m \cup H_g \cup H_b$ where $H_m$ comprises those copies that lie entirely inside some $V_i$, $H_g$ comprises those copies that intersect  every $V_i$ whose edge coloring coincides with the natural one given by the vertex partition (meaning the map from $R$ to $H$ takes vertex $i$ to a vertex in $V_i$), and $H_b$ comprises all other copies of $R$ (these include transversal copies but some vertex in any such copy will be in an inappropriate $V_i$). Let $h_m=|H_m|, h_g=|H_g|$ and $h_b=|H_b|$ so that
$$i(R; H) = h_m+h_g+h_b.$$
We will bound each of these three terms separately. First, note that since $n_i<n$, by induction
\begin{equation} \label{hmbound}
h_m \le \sum_j i(R; H[V_j]) \le \sum_j  g_k(n_j).
\end{equation}
Next we turn to $h_g$.
 Let $\Delta$ denote the number of $k$-sets that intersect each $V_i$ but are not counted by $h_g$. So a $k$-set counted by $\Delta$  either does not form a copy of $R$, or forms a copy of $R$ but its edge coloring does not coincide with the natural one given by the vertex partition $V_1 \cup \ldots \cup
 V_k$.
  Then
 \begin{equation} \label{hgbound}
 h_g =\prod_i n_i -\Delta
\end{equation}
and we need to bound $\Delta$ from below.
Note that the color or orientation of some pair in every member of $\Delta$ does not align with the implicit one given by our partition. With this in mind, let $D_{ij}$ be the set of pairs of vertices $\{v, w\}$ where $v \in V_i, w \in V_j$
such that either the color or orientation of $vw$ does not match that of $ij$ in $R$.
Let $\delta_{ij} =|D_{ij}|/{n \choose 2}$,  $D = \cup_{ij} D_{ij}$ and
$\delta= |D|/{n \choose 2}$. Let us  lower bound $\Delta$ by counting the misaligned
pairs from $D$ and then choosing the remaining $k-2$ vertices, one from each of the remaining parts $V_{\ell}$. This gives, for each $i<j$,
$$\Delta \ge |D_{ij}| \prod_{\ell \ne i,j} n_{\ell}=  \delta_{ij} {n \choose 2} \prod_{\ell \ne ij} n_{\ell}= \delta_{ij} {n \choose 2} \frac{\prod_{\ell} n_{\ell}}{n_in_j}.$$
Since $\sum_{ij} \delta_{ij} {n \choose 2} = \sum_{ij}|D_{ij}|=|D|=\delta{n \choose 2},$
we obtain by summing over $i,j$,
$$\Delta\left(\sum_{1\le i<j\le k} n_in_j\right) \ge \delta {n \choose 2} \prod_{\ell} n_{\ell}$$ which gives
\begin{equation} \label{Delta+}
h_g \le\prod_{\ell} n_{\ell}\left(1-\frac{\delta{n \choose 2}}{\sum_{1\le i<j\le k}n_in_j}\right) = \prod_{\ell} n_{\ell}\left(1-\frac{\delta{n \choose 2}}{{n \choose 2} - \sum_i {n_i \choose 2}}\right).
\end{equation}
Our next task is to upper bound $h_b$. For a vertex $x$ and $j \in [k]$, recall that $N_j(x) \subset V(H)$ is the set of  $y$ such that $x,y$ lie in a copy of $R$ with $x$ playing the role of vertex $j$ in $R$.  For $x \in V_i$, let $$Z(x) \df \max_{j \ne i} |N_j(x)| \qquad \hbox{ and } \qquad  z_{k,n} \df \max_{x \in V(H)} \frac{Z(x)}{(n-1)}.$$
Later we will give upper bounds for $z_{k,n}$. For now,
let us enumerate the set $J$ of tuples $(v,w,f)$ where $e=\{v,w\}\in D, f \in H_b$, $e \subset f$, and  $v \in V_i$, but $v$
plays the role of vertex $i'\ne i$ in the copy $f$ of $R$. In particular, all $k-1$ pairs $(v,x)$ with $x \in f$ contain color
$i'$. For $m=(v,w, f)\in J$, say that $m$ is 2-sided if $(w,v,f)\in J$ as well; otherwise say that $m$ is 1-sided.  Let $J_i$ be the set of $i$-sided tuples ($i=1,2$). We consider the weighted sum
$$S = 2|J_1| + |J_2|.$$
Observe that each $f \in H_b$ contains at least $k-2$ pairs from $D$. Indeed, if $f$ is transversal, then it must contain a  miscolored vertex which yields at least $k-1$ pairs from $D$ in $f$. If $f$ is not transversal, then take a largest color class $C$ of $f$, observe that at least $|C|-1$ of the vertices in $C$ are miscolored, and this yields at least $(|C|-1)(k-|C|) \ge k-2$ pairs from $D$ in $f$.
We conclude that each $f \in H_b$  contributes at least $2(k-2)$ to $S$ since $f$ contains at least $k-2$ pairs $e=\{u,v\}\in D$
and if $(v,w,f)$ is 1-sided it contributes 2 to $S$ while if it is 2-sided then it contributes 2 again since both $(v,w,f)$
and $(w,v, f)$ are counted with coefficient 1. This yields
\begin{equation} \label{Slower} S \ge 2(k-2)h_b.\end{equation}
On the other hand, we can bound $S$ from above by first choosing $e\in D$ and then $f \in H_b$ as follows. Call $v \in e =\{v,w\} \in D$ {\em correct in} $e$ if $v \in V_i$, $vw$ has color $\{i,j\}$ for some $j$ and  $v \to w$ in $H$ iff $i \to j$ in $R$; if $v$ is not correct in $e$ then say that $v$ is {\em wrong  in} $e$. The definition of $D$ implies that every $e\in D$ has at least one wrong vertex in $e$ (and possibly two wrong vertices). Let
$$D_i = \{\{v, w\}\in D: \{v, w\}\ \hbox{contains exactly $i$ wrong vertices}\} \qquad \hbox{($i=1, 2$)}.$$
The crucial observation is that
\begin{equation} \label{isided}( v,w,f) \in J_i \qquad \Longrightarrow \qquad \{v,w\} \in D_i \qquad \qquad \hbox{($i=1, 2$)}.\end{equation}
The reason this holds is that the color and orientation of $e \in D$ completely determine
the role that  its endpoints play in every copy of $R$ containing $e$.

Now, to bound $S$ from above, we use (\ref{isided}) and start by choosing $\{v,w\} \in D_i$ with wrong vertex $v$ and then the remaining $k-2$ vertices of $f\setminus e$. If $v \in V_i$, then, since $v$ is wrong in $e$ and $e \subset f \in H_b$,  the remaining $k-2$ vertices of $f \setminus e$ must all lie in $N_j(v)\setminus\{w\}$ for some $j \ne i$. So the number of choices for $f \setminus e$ is at most
$$p(|N_j(v)|-1, k-2) \le p(Z(v), k-2) \le p((n-1)z_{k,n}, k-2)$$ and for each choice of $f\setminus e$, we obtain  $m=(v,w,f) \in J_i$. This gives
$$ S \le 2 \sum_{\{v,w\} \in D_1}p((n-1)z_{k,n}, k-2) + 2\sum_{\{v,w\} \in D_2}p((n-1)z_{k,n}, k-2)=
2 \, |D| \, p((n-1)z_{k,n}, k-2)
.$$
Continuing, we obtain
\begin{equation} \label{Supper}
S \le 2 \, |D| \, p((n-1)z_{k,n}, k-2)\le 2 \, \delta{n \choose 2}\left(\frac{z_{k,n}}{k-2}\right)^{k-2} (n-1)^{k-2}.\end{equation}
Finally, (\ref{Slower}) and (\ref{Supper}) give
\begin{equation} \label{Delta-}h_b \le \frac{S}{2(k-2)} < \frac{ \delta{n \choose 2}}{k-2}\left(\frac{z_{k,n}}{k-2}\right)^{k-2}(n-1)^{k-2},\end{equation}
which is a refined version of \eqref{eq:bound_on_rb}.
Using (\ref{hmbound}), (\ref{Delta+}) and (\ref{Delta-})  we now have
$$i(R; H) < \sum_{i} g_k(n_i) +
\prod_{\ell} n_{\ell}\left(1-\frac{\delta{n \choose 2}}{{n \choose 2} - \sum_i {n_i \choose 2}}\right)
 + \frac{ \delta{n \choose 2}}{k-2}\left(\frac{z_{k,n}}{k-2}\right)^{k-2}(n-1)^{k-2}.$$
  Our final task now is to upper bound the RHS  by $g_k(n)$.

Since $\delta{n \choose 2} \le \sum_{ij}n_in_j = {n \choose 2} - \sum_i {n_i \choose 2}$, we have $\delta \in I\df (0, 1 - \sum_i {n_i \choose 2}/{n \choose 2})$. Viewing the expression above as a linear function of $\delta$, it suffices to check the endpoints of $I$.

If we let $\delta=0$, then
$$i(R; H) \le \sum_{i} g_k(n_i) + \prod_i n_i \le g_k(n)$$
and we are done. The last inequality holds since $g_k(n)$ is the maximum over all partitions of $n$, possibly with empty parts (as long as no part has size $n$), and $n_1+\cdots +n_k=n$ is one such partition.

If $\delta=1 - \sum_i {n_i \choose 2}/{n \choose 2}$ then we get
\begin{equation} \label{final} i(R; H) \le \sum_{i} g_k(n_i) +\frac{\sum n_in_j}{k-2}\left(\frac{z_{k,n}}{k-2}\right)^{k-2}(n-1)^{k-2}.\end{equation}

In order to show that the RHS in (\ref{final})  is at most $g_k(n)$, we will use explicit upper bounds on $g_k(n_i)$ and $z_{k,n}$ and lower bounds on $g_k(n)$.
Our first step is to state the following nontrivial lower bounds for $p(n,k)$. A proof is presented in the Appendix.

\medskip

\begin{lemma} \label{pbounds}
	For integers $k \ge 3$ and $n>k(k-1)$,
	$$ \left( \frac{n}{k}\right)^{k} \left( 1- e_k(n)\right)\le p(n,k)\le \left( \frac{n}{k}\right)^{k},$$
	where
	$e_k(n) = (4/27) (k^3/n^2)$.
\end{lemma}

We now give a bound on $z_{k,n}$.
\medskip

\begin{lemma} \label{zkn}
	For $k \ge 4$, $n>k(k-1)$, and $m=n-\lceil n/k \rceil$, let $z_{k,n}'$ be the largest real number $z \in (0, 1/2)$ that satisfies
	\begin{equation} \label{zeq} 	z^{k-1} + (1-z)^{k-1} \ge
	\frac{(k-1)^{k-1}}{k^{k-1}} \left(1-\frac{(4/27)(k-1)^3}{m^2}\right).\end{equation}
	
	Then	  $z_{k,n} \le z_{k,n}'$. Furthermore, $z_{k,n} <0.2611$ if either $k=4,\ n\geq 100$ or $k\geq 5$.
\end{lemma}

\proof
We begin by recalling  that $d_{min} >  g_k(n)-g_k(n-1)\ge p(m, k-1)$ where the second inequality holds by Lemma~\ref{trivial}. Recall  that $z_{k,n}=\max_y Z(y)/(n-1)$ and let $i \in [k]$ so that $x \in V_i$ achieves this maximum. Then by (\ref{eq_dn}) we have
\begin{equation}\label{eq_dmin} p(m, k-1) < d_{min} \le d(x) = \sum_{\ell=1}^k d_\ell(x) \le \sum_{\ell=1}^k p(|N_{\ell}(x)|, k-1).\end{equation}

Let $j \ne i$ be such that  $z_{k,n}=Z(x)/(n-1)=|N_{j}(x)|/(n-1)$. Then, writing $z=z_{k,n}$,  (\ref{eq_dmin}) continues as
\begin{eqnarray}
\nonumber  p(m, k-1) &<& \sum_{\ell=1}^k d_\ell(x) \le \sum_{\ell=1}^k p(|N_{\ell}(x)|, k-1)\leq p(|N_j(x)|,k-1)+p(n-1-|N_j(x)|,k-1)\\
\label{dmin}  & =& p(z(n-1), k-1) + p((1-z)(n-1), k-1).
\end{eqnarray}
Using Lemma~\ref{pbounds} this gives
$$\left( \frac{m}{k-1}\right)^{k-1} \left( 1- e_{k-1}(m)\right)
<\left( \frac{z(n-1)}{k-1}\right)^{k-1}
+\left( \frac{(1-z)(n-1)}{k-1}\right)^{k-1}.$$
Since $\lceil n/k \rceil \le (n+k-1)/k$, we have $m/(n-1) \ge (k-1)/k$
and the expression for $e_{k-1}(m)$ in Lemma~\ref{pbounds}  gives (\ref{zeq}).

The RHS of (\ref{zeq}) increases with $n$ and it is easy to see that it is $>2^{2-k}$ (the value of the LHS at $z=1/2$) already when $n=k(k-1)+1$. Hence the corresponding equation has two roots $0<z'<1/2< z''<1$ in the interval $(0,1)$ and $z_{k,n} \not\in I_{k,n} = (z', z'')$ which is an interval symmetric around 1/2. Since $I_{k, n+1} \supset I_{k,n}$, we conclude (for $k\geq 5$) that $z_{k,n} \not\in I_{k, k(k-1)+1}$. Direct calculation shows that $I_{4,100}\supset (0.2611,0.7389),\ I_{5,21}\supset (0.2611,0.7389)$, and it is an easy exercise to see that the intervals $I_{k,k(k-1)+1}$ only grow with $k$.
To complete the proof, we only need to show that $z_{k,n} \le 1/2$ for $k \ge 4$ and $n>k(k-1)$.

Suppose for the sake of contradiction that  $Z(x)= |N_j(x)|\ge (n-1)/2$ (recall that $x \in V_i$). Then, since $x \in V_i$, we have
$$
d_j(x)\leq d_i(x) \leq p((1-z)(n-1), k-1)
$$
and hence instead of the bound $d_j(x)\leq p(z(n-1),k-1)$ we could use in \eqref{dmin} this better bound. That would give us
$$
p(m,k-1) < 2p((1-z)(n-1),k-1) \leq 2p(\lfloor(n-1)/2\rfloor,k-1).
$$
This, however is false e.g. since, as we argued above, $1/2\in I_{k,n}$.
This contradiction shows that in fact $z_{k,n}\leq 1/2$ and completes the proof of Lemma \ref{zkn} \qed
\medskip

 Our next lemma provides  bounds for $g_k(n)$. We will give a proof in the Appendix.
 \medskip

\begin{lemma} \label{gbounds}
	For $k \ge 4$ and $n>k(k-1)$
	
	$$\frac{n^k - k^3 n^{k-2}}{k^k-k} \le  g_k(n) \le \frac{n^k-n}{k^k-k}.$$
\end{lemma}

Our final task is to provide a nontrivial upper bound on each $n_i$. Write $p_i=n_i/n$ and  $e'_k(n)= k^3/n^2$.  Assume w.l.o.g. that $p_1=\max p_i$.
\medskip

\begin{lemma} \label{pi} Let $k \ge 4$ and $n>\begin{cases}k(k-1) & \text{if}\ k\geq 5\\ 100 & \text{if}\ k=4\end{cases}$. Then we can assume w.l.o.g. that
	$p_1<0.86$.
	\end{lemma}

\proof  Let $p_1=1-c_k$ and assume w.l.o.g. that $p_k=\min_{i>1} p_i$.  Then our assumption implies that  $p_k \le c_k/(k-1)$. We consider
\begin{equation} \label{irh} i(R; H) = \sum_{v \in V} d_k(v) = \sum_{v \not\in V_k} d_k(v) + \sum_{v \in V_k} d_k(v).\end{equation}
Note that for $v \not\in V_k$, $|N_k(v)| \le z_{k,n}(n-1)$ and so $d_k(v)\le p(z_{k,n}(n-1), k-1)$. For $v \in V_k$, we will use the weaker bound $d_k(v)\le p(n-1, k-1)$. As  we may assume that $i(R; H) \ge g_k(n)$ (otherwise we are done by induction), Lemma~\ref{gbounds} and (\ref{irh}) give
\begin{equation} \label{save} \frac{n^k(1-e'_k(n))}{k^k-k} \le g_k(n) \le i(R; H) \leq n \, p(z_{k,n}(n-1), k-1) + \frac{c_k n}{k-1} p(n-1, k-1).\end{equation}
Dividing by $n^k$ and using $p(n-1, k-1) < p(n,k-1) \le (n/(k-1))^{k-1}$ we obtain
$$c_k \ge \frac{(k-1)^k}{k^k-k}(1-e'_k(n)) - (k-1)z_{k,n}^{k-1} \ge
\frac{(k-1)^k}{k^k-k}\left(1-\frac{k^3}{(k(k-1)+1)^2}\right) - (k-1)z_{k,n}^{k-1}
.$$
By the last part of Lemma \ref{zkn}$, z_{k,n}\leq 0.27$. This shows that $c_4>0.14$ or $p_1=1-c_4<0.86$, and it is an easy matter to see that the bound only improves as $k$ increases.
	\qed

\medskip

We are now ready to complete the proof.
Recall that our main equation is
\begin{equation} \label{final2} i(R; H) \le \sum_{i} g_k(n_i) +\frac{\sum n_in_j}{k-2}\left(\frac{z_{k,n}}{k-2}\right)^{k-2}(n-1)^{k-2},\end{equation}
and we are to show that the RHS is at  most $g_k(n)$.

We treat the case $k=4,\ n\leq 100$ by exhaustive search through all partitions. A simple Maple worksheet verifying this fact (as well as a few other numerical facts that we state in our proof) can be found at the web pages   {\tt http://homepages.math.uic.edu/\~{}mubayi/papers/ErdosHajnalmw.pdf} and {\tt http://people.cs.uchicago.edu/\~{}razborov/files/ErdosHajnal.mw}. Thus,
in what follows we always assume that $k=4$ entails $n\geq 100$. In particular, we can utilize the conclusions of Lemmas \ref{zkn} and \ref{pi}:
$$
z_{k,n} < 0.2611,\ p_1<0.86.
$$

Dividing by $n^k/(k^k-k)$ and using Lemma~\ref{gbounds}, we see that it suffices to prove
$$L \df \sum_i p_i^k + A \sum_{ij} p_ip_j  \le 1-e'_k(n)$$
where as before $e'_k(n) = k^3/n^2$ and
$$
A =  A(k) \df \frac{(k^k-k)\cdot 0.2611^{k-2}}{(k-2)^{k-1}}.
$$
Since $p_1\ge p_i$ for all $i$ and $\sum_i p_i = 1$,
$$L \le p_1^{k-2} \left(\sum_i p_i^2\right) + A \sum_{ij}p_ip_j = p_1^{k-2}\left(1-2\sum_{ij}p_ip_j\right) + A \sum_{ij}p_ip_j = p_1^{k-2} + (A-2p_1^{k-2})\sum_{ij} p_ip_j.$$

 If $A< 2p_1^{k-2}$, then $k>4$ since $A(4)>2.14$ and the coefficient of $\sum_{ij}p_ip_j$ is negative. Lemma~\ref{pi} then gives
\begin{equation} \label{L}
L < p_1^{k-2} \le p_1^3 <(0.86)^3  < 0.7 < 1-e'_4(21)\le 1-e'_k(n).
\end{equation}
Now assume that $A \ge 2p_1^{k-2}$.  Then, using $\sum_{ij}p_ip_j\le (k-1)/2k$ it is enough to show
\begin{equation} \label{Acalc}
p_1^{k-2} + \frac{k-1}{2k}(A-2p_1^{k-2}) < 1-e'_k(n).
\end{equation}
For the same reasons as at the end of the proof of Theorem \ref{asymptotic}, we must split further analysis into two cases.

If $k\geq 5$, we apply the trivial bound $\frac{k-1}{2k}<\frac 12$ that reduces \eqref{Acalc} to merely
$$
A <2 (1-e_k'(n)).
$$
This holds since
$$A(k) \le A(5) < 1 < 2(1 - e'_5(21))< 2(1 - e'_k(n)).$$

For $k=4$, \eqref{Acalc} becomes
$$
\frac{p_1^2}{4}+\frac 38A(4) <1-e'_4(n).
$$
In this case $p_1<0.86$, $A(4)<2.15$ and
$$
\frac{p_1^2}{4}+\frac 38A(4) <0.992 <1-e_4'(100).
$$

The proof of Theorem~\ref{main} is complete.
 \qed
\bigskip

\section{Appendix}

Here we give the proofs of Lemmas~\ref{trivial}, \ref{pbounds} and \ref{gbounds}.

\subsection{Proof of Lemma~\ref{trivial}}
We only consider partitions $n_1+\cdots+ n_k=n$ where $0\le n_i<n$ for all $i$; the partition is equitable if $|n_i-n_j|\le 1$ for all $i\ne j$.
For a vertex $v$ in a hypergraph $H$, we write $d_H(v)$ for the degree of $v$ in $H$. We will use the notation  $H$ for the edge set of $H$.

\begin{definition}
	For $k \ge 3$ and $n \ge 0$, let $\gkn$ be the family of $n$-vertex $k$-graphs defined inductively as follows: For $n<k$, $\gkn$ comprises the single $n$-vertex $k$-graph with no edge. For $n\ge k$, the vertex set $V$ of any $G \in \gkn$  is partitioned as $V_1 \cup \ldots \cup V_k$ with $n_i:=|V_i|$ and $0\le n_i<n$ for all $i$. For the edge set,  $G[V_i] \in \gkni$ for each $i$, and in addition $G$ contains all edges that have one point in each $V_i$. Call  $V_1 \cup \ldots \cup V_k$ the {\em defining partition} for $G$, or simply, {\it the partition} for $G$.
\end{definition}

\begin{definition}
	Let $H_k(n) \in \gkn$ be the following $k$-graph. For $n<k$, $H_k(n)$ is the unique member of $\gkn$, and for $n\ge k$, the defining partition  $V_1 \cup \ldots \cup V_k$ of  $H_k(n)$ is an equitable partition ($||V_i|-|V_j||\le 1$ for all $i \ne j$) and for each $i$, the subgraph induced by $V_i$ is isomorphic to $H_k(|V_i|)$.  We let $h_k(n):=|H_k(n)|$.
\end{definition}

Our proof of Lemma~\ref{trivial} will use induction on $n$ and so we need one more definition.

\begin{definition}
	A vertex $v$ of $G \in \gkn$ is {\em G-good} if the following holds: for $n<k$ every vertex is $G$-good. For $n\ge k$, if  $V_1 \cup \ldots \cup V_k$ is the partition for $G$, and $|V_i|\ge |V_j|$ for all $j$, then $v$ is $G$-good  if $v \in V_i$ and $v$ is  $G[V_i]$-good.   In other words, a vertex is $G$-good if it lies in a largest part $V_i$ in the partition for $G$ and the same is true inductively within $V_i$.
\end{definition}
Removing any vertex $v$ from $G \in \gkn$ results in a $k$-graph $G-v \in {\cal G}_k(n-1)$.
Moreover,  if $v$ is $H_k(n)$-good, then
\begin{equation}\label{hkn} H_k(n)-v \cong H_k(n-1).\end{equation}
Indeed, if we remove $v$ from $H_k(n)$, then the  partition  for $H_k(n)$, after removal of $v$,  is still equitable and the same remains true of all inductively defined partitions.  Now (\ref{hkn}) shows that every two  $H_k(n)$-good vertices have the same degree and hence we may define $\delta_k(n) = d_{H_k(n)}(v)$ where $v$ is any  $H_k(n)$-good vertex. Observe that
\begin{equation} \label{star}
\delta_k(n) = d_{H_k(n)}(v) = d_{H_k(\lceil n/k \rceil)}(v) + p(n-\lceil n/k\rceil, k-1) = \delta_k(\lceil n/k \rceil) +p(n-\lceil n/k\rceil, k-1).
\end{equation}
Finally, (\ref{hkn}) gives
\begin{equation} \label{hkn2}h_k(n-1)+ \delta_k(n) = h_k(n).\end{equation}
\smallskip

\begin{lemma} \label{dlem}
	Let $G \in \gkn$ and $v$ be $G$-good. Then $d_G(v) \le \delta_k(n)$.
\end{lemma}

\proof
Proceed by induction on $n$. The cases $n<k$ are trivial since $d_G(v) =0= \delta_k(n)$. Let $V_1 \cup \ldots \cup V_k$ be the partition for $G$ and $n_i:=|V_i|$ with $n_1\ge n_2\ge \cdots \ge n_k$ and assume wlog that $v \in V_1$.
Let $X_1 \cup \cdots  \cup  X_k$ be the partition for $G[V_1] \in {\cal G}_k(n_1)$, $x_i:=|X_i|$ with $b:=x_1\ge \cdots  \ge x_k$  and assume wlog that $v \in X_1$. Note that $b \ge \lceil n_1/k \rceil$.
Let
$$a_1:= |V_2 \cup \cdots \cup V_k|=n-n_1, \qquad a_2:= |V_1|-|X_1| = n_1-b.$$
Since $v$ is $G[X_1]$-good, $G[X_1] \in {\cal G}_k(b)$, and $b<n$,  induction implies
$d_{G[X_1]}(v) \le \delta_k(b)$ and hence
\begin{equation} \label{2star}
d_G(v) = d_{G[V_1]}(v) + \prod_{j=2}^k n_j = d_{G[X_1]}(v)+\prod_{\ell=2}^k x_{\ell} +\prod_{j=2}^k n_j \le
\delta_k(b) + p(a_1, k-1) + p(a_2, k-1).\end{equation}

{\bf Case 1. $b  \le  \lceil n/k \rceil$.}
Note that $a_i +b \ge \lceil n/k \rceil$ for $i=1,2$ since $a_1+b=(n-n_1)+b\ge (n-n_1)+n_1/k\ge n/k$ and $a_2+b = n_1\ge n/k $. For fixed $b$ and $n$,  $a_1+a_2=n-b$ is also fixed, so by Lemma~\ref{basics}c),
$p(a_1, k-1)+ p(a_2,k-1)$ is uniquely maximized when $a_1$ or $a_2$ is as small possible, namely  $\{a_1, a_2\}=\{\lceil n/k \rceil-b,n-\lceil n/k \rceil\}$ where we use the assumption $\lceil n/k \rceil-b\ge 0$. Consequently,
\begin{equation} \label{d} d_G(v)  \le
\delta_k(b) +  p(\lceil n/k \rceil-b, k-1)+p(n-\lceil n/k \rceil, k-1).\end{equation}
If $b=\lceil n/k \rceil$, then (\ref{star}) and (\ref{d}) give $d_G(v) \le \delta_k(n)$ and we are done, so assume that $b <  \lceil n/k \rceil$.
Consider $K \in {\cal G}_k(\lceil n/k \rceil)$ whose defining partition has largest part $B$ of size $b$ and all other $k-1$ parts form an equitable partition of $\lceil n/k \rceil -b$.
Since $b =x_1 \ge \lceil n_1/k\rceil  \ge \lceil \lceil n/k \rceil/k\rceil$,  $B$ is indeed a largest part. Further, let $K[B]\cong H_k(b)$ (for all other parts $C$, choose $K[C]$ arbitrarily) and let $w$ be a $K[B]$-good vertex. Then $d_{K[B]}(w)=\delta_k(b)$ and $w$ is also $K$-good, so by induction,
$$\delta_k(b) + p(\lceil n/k \rceil-b, k-1) = d_{K[B]}(w) + p(\lceil n/k \rceil-b, k-1) = d_{K}(w) \le \delta_k(\lceil n/k \rceil).$$
Continuing (\ref{d}) we obtain
$$d_G(v) \le \delta_k(\lceil n/k \rceil) + p(n-\lceil n/k \rceil, k-1) = \delta_k(n)$$
where we use (\ref{star}) for the last equality.

{\bf Case 2. $b > \lceil n/k \rceil$.}  In this case (\ref{2star}) and Lemma~\ref{basics} yield
$$d_G(v)  \le
\delta_k(b) + p(a_1, k-1) + p(a_2, k-1) \le \delta_k(b) +p(a_1+a_2, k-1) =
\delta_k(b) +p(n-b, k-1).$$
Consider $G' \in {\cal G}_k(n)$ whose defining partition has largest part $B'$ of size $b$ and all other $k-1$ parts form an equitable partition of $ n-b$.
Since $b \ge \lceil n/k\rceil $,  $B'$ is indeed the largest part. Further, let $G'[B']\cong H_k(b)$ and let $X_1' \cup \ldots\cup  X_k'$ be the partition for  $G'[B']$ with $|X_1'| =\lceil b/k \rceil$.
Let $v'\in X_1'$ be a $G'[B']$-good vertex. Since $B'$ is the largest part in the partition for $G'$,  $v'$ is also $G'$-good. Since $|X_1'| = \lceil b/k \rceil  \le  \lceil n/k \rceil$,  by the proof in Case 1 (with $(G', X_1', v')$ playing the role of $(G, X_1, v)$) we conclude that $d_{G'}(v') \le \delta_k(n)$. On the other hand, since $v'$ is $G'[B']$-good and $G'[B']\cong H_k(b)$,
$$d_{G'}(v')=d_{G'[B']}(v') + p(n-b, k-1) = d_{H_k(b)}(v') + p(n-b, k-1) = \delta_k(b)
+p(n-b, k-1).$$
We therefore have
$$d_G(v) \le \delta_k(b) + p(n-b, k-1) = d_{G'}(v') \le \delta_k(n). \qed $$

{\bf Proof of Lemma~\ref{trivial}.}  We are to show that $g_k(n)=h_k(n)$. In other words, for each $G \in \gkn$ we must show that $|G| \le h_k(n)$. We proceed by induction on $n$. The cases $n\le k$ are trivial, so assume $n > k$.  Pick a $G$-good vertex $v$. Then $G-v \in {\cal G}_k(n-1)$ and by induction, Lemma~\ref{dlem}, and (\ref{hkn2}),
$$|G| =|G-v| + d_G(v) \le h_k(n-1) + \delta_k(n) = h_k(n). \qed$$

\subsection{Proof of Lemma~\ref{pbounds}}

 We begin with the following inequality.

\begin{lemma} \label{ablemma}
	Let $a,b,k,n$ be positive integers with
	$k =a+b \ge 2$, and $n \ge \max\{a,b\}$. Then
	$$\left(1+\frac{a}{n}\right)^b \left(1-\frac{b}{n}\right)^a \ge 1 - \frac{\max\{ab^2, ba^2\}}{n^2} \ge 1-  \frac{(4/27)k^3}{n^2}.$$
\end{lemma}
\proof Suppose first that $a\le b$. Consider $a$ copies of the number $1+b/n$ and $b-a$ copies of the number 1. The arithmetic mean of these numbers is $1+a/n$ hence
the AMGM inequality gives
$$\left(1+\frac{a}{n}\right)^b \ge \left(1+\frac{b}{n}\right)^a \cdot 1^{b-a}$$
and Bernoulli's estimate yields
$$\left(1+\frac{a}{n}\right)^b  \left(1-\frac{b}{n}\right)^a
\ge \left(1-\frac{b^2}{n^2}\right)^a
\ge 1-\frac{ab^2}{n^2}.$$
If $ a \ge b$, then a similar argument applies by taking $b$ copies of $1-a/n$ and $a-b$ copies of 1.
\qed

\bigskip

{\bf Proof of Lemma~\ref{pbounds}.} We are to show that
for $k \ge 3$ and $n>k(k-1)$,
$$\left( \frac{n}{k}\right)^{k} \left( 1- e_k(n)\right)\leq p(n,k)\leq \left( \frac{n}{k}\right)^{k},$$
where
$e_k(n) = (4/27) (k^3/n^2)$.
The upper bound for $p(n,k)$ is trivial so we only prove the lower bound. Let $n \equiv t$ (mod $k$) where $0 \le t <k$. Then
$$p(n,k) = \left( \frac{n}{k} + \frac{k-t}{k} \right) ^t \left( \frac{n}{k} - \frac{t}{k} \right)^{k-t}=\left(\frac{n}{k}\right)^k\left(1 + \frac{k-t}{n} \right) ^t \left( 1 - \frac{t}{n} \right)^{k-t}.$$ Now apply Lemma~\ref{ablemma} with $a=k-t$ and $b=t$.  \qed

\subsection{Proof of Lemma~\ref{gbounds}}

We are to show that for $k \ge 4$ and $n>k(k-1)$
$$\frac{n^k - k^3 n^{k-2}}{k^k-k} \le  g_k(n) \le \frac{n^k-n}{k^k-k} .$$
The upper bound is by induction on $n$ and we prove it for all $n\ge 1$. The base cases $n \le k$ are obvious so let $n>k$. For the induction step, apply Lemma~\ref{trivial} and take the equitable partition $n=\sum_i n_i$ that achieves the definition of $g_k(n)$. Then each $n_i<n$ and by induction,
$$g_k(n) = \sum_i g_k(n_i) + \prod_i n_i \le \sum_i \frac{n_i^k-n_i}{k^k-k} + \prod_i n_i
=\sum_i \frac{n_i^k}{k^k-k} -\frac{n}{k^k-k}+ \prod_i n_i \le \frac{n^k-n}{k^k-k}$$
where the last inequality  (after dividing by $n^k/(k^k-k)$) is   (\ref{eq:left_end}).

For the lower bound, we take  an equitable partition $\sum n_i=n$  and proceed by induction on $n$.  Let us first assume that $k \ge 5$. We will actually prove that
$$(k^k-k) \, g_k(n) \ge \begin{cases}
n^k-k^4 n^{k-2} \qquad \hbox{ if $n\le k(k-1)$}\\
n^k - k^3 n^{k-2}\qquad \hbox{ if $n > k(k-1)$.}
\end{cases}
$$
Note that the first bound is trivial for $n \le k(k-1)$ since it is negative.  For $n>k(k-1)$,  Lemma~\ref{pbounds} and induction imply that   $g_k(n)$ is at least
$$ \sum_i g_k(n_i) + p(n,k)\ge \sum_i\frac{n_i^k-k^4n_i^{k-2}}{k^k-k}+\left(\frac{n}{k}\right)^k\left(1-e_k(n)\right)\ge \frac{n^k}{k^k-k} -\left(
\frac{\sum k^4n_i^{k-2}}{k^k-k} +\frac{n^k e_k(n)}{k^k}\right)$$
where we used Jensen's inequality to obtain
$$\sum \frac{n_i^k}{k^k-k} + \left(\frac{n}{k}\right)^k \ge \left(\frac{n}{k}\right)^k\left(\frac{k}{k^k-k} + 1\right) =
\frac{n^k}{k^k-k}.$$
Since $n>k(k-1)$ and the $n_i$'s are an equitable partition,
\begin{equation} \label{k4} \sum_{i=1}^k n_i^{k-2} < \sum_{i=1}^k (n/k+1)^{k-2} = k(n/k+1)^{k-2} = \frac{n^{k-2}}{k^{k-3}}(1+k/n)^{k-2} <\frac{n^{k-2}}{k^{k-3}}\cdot e.\end{equation}  Using $k \ge 5$, we now  have
\begin{equation} \label{1827}
\frac{\sum k^4n_i^{k-2}}{k^k-k} \le e \cdot \frac{k^4}{k^k-k} \cdot \frac{n^{k-2}}{k^{k-3}} < \frac{18}{27} \frac{k^3 n^{k-2}}{k^k-k}.\end{equation}
Hence

\begin{equation} \label{finalapp}\frac{\sum k^4n_i^{k-2}}{k^k-k} +\frac{ n^k e_k(n)}{k^k} <
\frac{ (18/27)k^3 n^{k-2}}{k^k-k}  +\frac{ n^k (4/27)(k^3/n^2)}{k^k-k} <
\frac{k^3n^{k-2}}{k^k-k}\end{equation}
and the proof is complete.

 Now we  assume that $k=4$. In this case we show by induction on $n$ that for all $n>0$, $$g_4(n) \ge \frac{n^k-k^3n^{k-2}}{k^k-k} = \frac{n^4-64n^2}{252}.$$ The cases $n \le 8$ are trivial since the RHS is nonpositive so assume that $n \ge 9$. The cases $9\le n \le 12$ can be checked by direct computation, so assume that $n>12$ and we can apply the bounds in Lemma~\ref{pbounds}. We proceed as in the proof for $k \ge 5$ except that all occurrences of the  $k^4$ term are replaced by $k^3=4^3$. Now (\ref{k4}) becomes
	$$\sum_{i=1}^4 n_i^2 < \frac{n^2}{4}(1+4/n)^2\le \frac{n^2}{4}(1+4/9)^2< (2.1) \frac{n^2}{4}.$$
	and, since $2.1<72/27$,  (\ref{1827}) becomes
	$$\frac{\sum 4^3 n_i^2}{252} \le  \frac{4^3}{252} \cdot (2.1) \cdot \frac{n^2}{4} < \frac{18}{27} \frac{4^3 n^2}{252}.$$
	Finally, (\ref{finalapp}) becomes
	$$\frac{\sum 4^3n_i^2}{252}
	+\frac{ n^4 e_4(n)}{4^4} <
	\frac{ (18/27)4^3 n^{2}}{252}  +\frac{ n^4 (8/27)(4^3/n^2)}{252} <
	\frac{64n^{2}}{252}$$
	and the proof is complete.\qed

\section*{Acknowledgment}

We are grateful to David Conlon and Sergey Tarasov for several useful remarks and to a referee for extremely careful and expert reading of the paper and providing a shorter proof of Lemma~\ref{ablemma}.

\end{document}